\documentclass{article}
\usepackage{amssymb,latexsym,amsmath,amsthm,mathrsfs}
\usepackage{graphicx}
\newcommand{\comment}[1]{}

\DeclareMathOperator{\tang}{tang}
\DeclareMathOperator{\Atang}{Atang}
\begin{document}
\title{The solution of a memorable problem by a special artifice of
calculation\footnote{Presented to the
St. Petersburg Academy on March 22, 1779.
Originally published as
{\em Solutio problematis ob singularia calculi artificia memorabilis},
M\'emoires de l'acad\'emie des sciences de
St-P\'etersbourg \textbf{2} (1810), 3--9.
E731 in the Enestr{\"o}m index.
Translated from the Latin by Jordan Bell,
Department of Mathematics, University of Toronto, Toronto, Canada.
Email: jordan.bell@utoronto.ca}}
\author{Leonhard Euler}
\date{}
\maketitle

1. The problem which I take up for solving here can be thus enunciated:
{\em To find a curved line $AM$} (Fig. 1), {\em with coordinates $CP=x,PM=y$,
the arc $AM=s$ and chord $CM=\sqrt{xx+yy}=z$, so that the integral
formula $\int vds$ obtains a maximum or a minimum value, with
$v$ any fixed function of $z$.}

2. In general, if a relation is sought between two variables
$x$ and $y$, and it is put $dy=pdx$, and $V$ is any function of
these $x,y$ and $p$, so that its differential has this form: $dV=Mdx+Ndy+Pdp$,
then the integral formula $\int Vdx$ will have a maximum or a minimum value
when it happens that $Ndx=dP$, so that this equation would express the sought
relation between $x$ and $y$.

3. While the entire problem is reduced to this equation, 
it will however be useful to further consider another equation which is
 equivalent
to the first.
For when $Ndx=dP$,
which multiplied by $p$ becomes $Ndy=pdP$,
 substituting in this value will produce
\[
dV=Mdx+Pdp+pdP=Mdx+d\cdot Pp.
\]
It therefore follows that
$Mdx=d\cdot (V-Pp)$. This equation is the most
accommodating for use in our analysis.

\begin{figure}
\begin{center}
\includegraphics{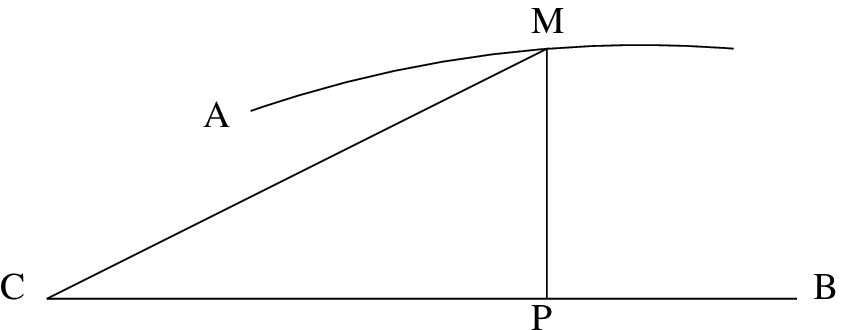}
\end{center}
\begin{center}
Fig. 1
\end{center}
\end{figure}

4. Let us now transfer these general precepts to the proposed problem.
And first putting $dy=pdx$ we will have $ds=dx\sqrt{1+pp}$.
Next, since it is
$z=\sqrt{xx+yy}$, it will be $dz=\frac{xdx+ydy}{z}$.
Then indeed, as $v$ is a function of $z$, let us put $dv=qdz$ and it will
then be $dv=\frac{q(xdx+ydy)}{z}$. 
Now therefore for the formula of the maximum or the minimum we will have
$V=v\sqrt{1+pp}$. By differentiating this we will get
\[
dV=\frac{q(xdx+ydy)\sqrt{1+pp}}{z}+\frac{vpdp}{\sqrt{1+pp}},
\]
and with the differential compared to the general form it will be
$M=\frac{qx\sqrt{1+pp}}{z}$, $N=\frac{qy\sqrt{1+pp}}{z}$ and
$P=\frac{vp}{\sqrt{1+pp}}$.

5. Therefore, the equation comprising the total solution
of our problem, which was $Ndx=dP$, now assumes this form:
$\frac{qydx\sqrt{1+pp}}{z}=d\cdot \frac{vp}{\sqrt{1+pp}}$.
As well, the other equation which will be used,
since $V-Pp=\frac{v}{\sqrt{1+pp}}$, may be expressed thus:
\[
\frac{qxdx\sqrt{1+pp}}{z}=d\cdot \frac{v}{\sqrt{1+pp}}.
\]
Hence, since the first equation is
\[
d\cdot \frac{vp}{\sqrt{1+pp}}=\frac{vdp}{\sqrt{1+pp}}+
p\cdot d\cdot \frac{v}{\sqrt{1+pp}}=
\frac{vdp}{\sqrt{1+pp}}+\frac{qxdy\sqrt{1+pp}}{z},
\]
we will be led to this:
\[
\frac{q(ydx-xdy)\sqrt{1+pp}}{z}=\frac{vdp}{\sqrt{1+pp}}
\]
and hence it will be
\[
ydx-xdy=\frac{vzdp}{q(1+pp)},
\]
and this is the equation which we employ for getting the solution of our
problem.

6. Let us divide this equation by $zz=xx+yy$, so that we will obtain
this form:
\[
\frac{ydx-xdy}{xx+yy}=\frac{vdp}{qz(1+pp)},
\]
where it is clear that the integral of the first side is $\Atang\frac{x}{y}$.
In truth it seems here that little gain be obtained, since the latter part
of the equation is completely intractable. In the meanwhile however let
us put $\varphi$ to be the angle whose tangent is $\frac{x}{y}$,
so that our equation would be $d\varphi=\frac{vdp}{qz(1+pp)}$.

7. With this angle $\varphi$ introduced, we will be able to render
the coordinates $x$ and $y$ susceptible to calculation. For since it is
$\frac{x}{y}=\Atang \varphi$, it will be $x=z\sin \varphi$ and $y=z\cos\varphi$,
by means of whose values the letter $p$ can also be extracted.
Indeed because $p=\frac{dy}{dx}$, it will be
\[
p=\frac{dz\cos\varphi-zd\varphi\sin\varphi}{dz\sin\varphi+zd\varphi\cos\varphi}.
\]
Now let us put $dz=tzd\varphi$, so that it becomes
\[
p=\frac{t\cos\varphi-\sin\varphi}{t\sin\varphi+\cos\varphi}=\frac{t-\tang\varphi}{1+t\tang\varphi}.
\]
This expression 
clearly expresses the tangent of the difference of two angles, the tangent of
the first of which is $=t$, while the latter angle is $=\varphi$.

8. Thus for $p$ equal to the tangent of any angle whatsoever, let us put
$p=\tang \omega$, and $\omega$ will be the difference of these
angles, namely $\omega=\Atang t-\varphi$, whence $d\omega=\frac{dt}{1+tt}-d\varphi$. 
Also indeed, because $p=\tang \omega$ and so $\omega=\Atang p$, it will further
be $d\omega=\frac{dp}{1+pp}$. Hence our equation to be resolved will
be $d\varphi=\frac{vd\omega}{qz}$. Moreover from the preceding,
since the form was $d\omega=\frac{dt}{1+tt}-d\varphi$, hence
\[
\frac{qzd\varphi}{v}=\frac{dt}{1+tt}-d\varphi \quad \textrm{or} \quad 
d\varphi(1+\frac{qz}{v})=\frac{dt}{1+tt}.
\]

9. Also, because we have put $dz=tzd\varphi$, it will be
$d\varphi=\frac{dz}{tz}$, and having substituted in this value our equation takes the form:
\[
\frac{dz}{z}(1+\frac{qz}{v})=\frac{tdt}{1+tt}.
\]
Then, since $qdz=dv$, the integration can be done most conveniently by
logarithms; for it will be
\[
lz+lv=l\sqrt{1+tt}-ln,
\]
and consequently we will obtain this integrated formula: $vz=\frac{\sqrt{1+tt}}{n}$.

10. Now let us investigate the value of $t$ from this equation, in which
it will be
$t=\sqrt{nnvvzz-1}$. For since $t=\frac{dz}{zd\varphi}$, we gather from
this equation that
\[
d\varphi=\frac{dz}{z\sqrt{nnvvzz-1}},
\]
which is a differential equation of the first degree between the angle
$\varphi$ and the distance $CM=z$, if indeed $v$ is a function of $z$ itself.
Indeed for the angle it is noted that $\tang\varphi=\frac{x}{y}$
and hence $x=z\sin\varphi$ and $y=z\cos\varphi$,
so that now the two coordinates $x$ and $y$ can by expressed by the same
variable $z$, which is the most complete solution of our problem.

11. Here I am compelled however to admit that I would perhaps not have
obtained this solution if it had not already been noted by me elsewhere;
and for this reason, the artifice which I used in this calculation seems worthy
of more attention, for it is not completely obvious and without doubt
will be able to be used in many other cases.

12. Thus I will treat here another solution of the same problem, which
I have employed to get the same final solution found above without
any detours. Namely, I recast the entire question into two other
coordinates, suited for defining a curve. The first is the distance $CM$ (Fig. 2), which I shall call here $=x$, and the other is the angle $BCM$,
designated by the letter $y$. Then an element of the curve will be
$Mm=ds=\sqrt{dx^2+xxdy^2}$, which by putting $dy=pdx$ turns into
$ds=dx\sqrt{1+ppxx}$, whence the integral formula for maximization
or minimization will be $\int vdx\sqrt{1+ppxx}$.

\begin{figure}
\begin{center}
\includegraphics{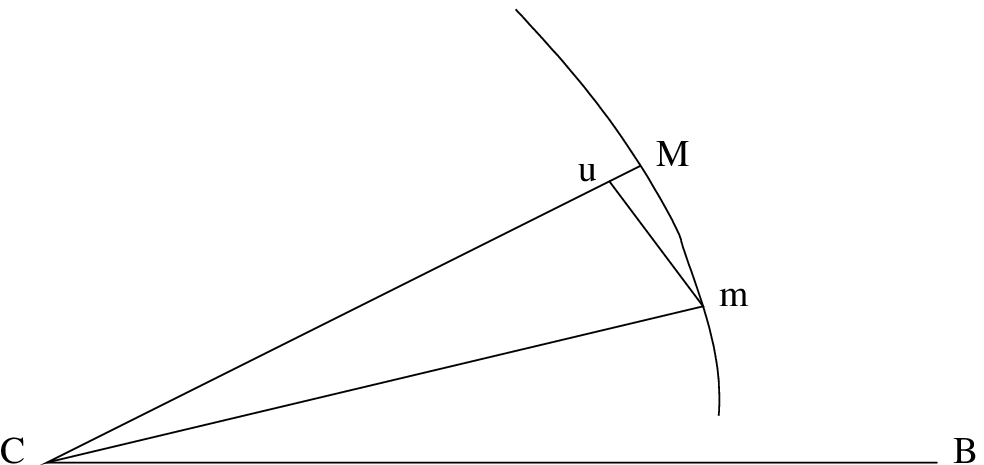}
\end{center}
\begin{center}
Fig. 2
\end{center}
\end{figure}

13. We then compare this formula with the general formula $\int Vdx$, and we will
have $V=v\sqrt{1+ppxx}$. Therefore this quantity, because $v$ is a function of $x$ only, involves only the variable quantities $x$ and $p$, with the third
$p$ completely excluded; whence, since we have put
$dV=Mdx+Ndy+Pdp$, it will be
\[
N=0 \quad \textrm{and} \quad P=\frac{vpxx}{\sqrt{1+ppxx}}.
\]
Hence the equation containing the solution, $Ndx=dP$, turns into
this:
$dP=0$,
whence
$P=\textrm{const.}=\frac{1}{n}$, so that
$nvpxx=\sqrt{1+ppxx}$ and hence it at once follows that
\[
p=\frac{1}{x\sqrt{nnvvxx-1}}=\frac{dy}{dx},
\]
and thus now we have arrived at this equation:
\[
dy=\frac{dx}{x\sqrt{nnvvxx-1}}.
\]

14. We may now transfer this very simple solution gained for the denominator
in the above solution, namely if we write $z$ in place of the letter $x$ and
the element $d\varphi$ in place of $dy$, and in this way the solution
of our problem will be contained in this equation:
\[
d\varphi=\frac{dz}{z\sqrt{nnvvzz-1}},
\]
which, since $v$ is a function of $z$ itself, completely agrees with that which
we deduced in the prior solution by many detours. In particular it should
be observed here that this solution is always valid, whatever function of $z$
is taken for $v$. In particular as well, this is noteworthy because if $v$ is taken to be a power of $z$, the satisfying curve will develop algebraically.

15. For let us put $v=z^\lambda$, and we will have this equation for the sought
curve:
\[
d\varphi=\frac{dz}{z\surd(nnz^{2\lambda+2}-1)}.
\]
For expanding this, let us set $\sqrt{nnz^{2\lambda+2}-1}=u$, so that it becomes
$d\varphi=\frac{dz}{zu}$. Then indeed it will be
\[
nnz^{2\lambda+2}=uu+1
\]
and by taking the logarithmic differentials
\[
(2\lambda+2)\frac{dz}{z}=\frac{2udu}{1+uu}
\]
and hence
\[
\frac{dz}{z}=\frac{udu}{(\lambda+1)(1+uu)},
\]
so that we will now have $(\lambda+1)d\varphi=\frac{du}{1+uu}$, and then
by integrating
\[
(\lambda+1)\varphi=A\tang u.
\]
But if then we take $\psi$ as the angle whose tangent is $\sqrt{nnz^{2\lambda+2}-1}$, it will be $(\lambda+1)\varphi=\psi$ and hence $\varphi=\frac{\psi}{\lambda+1}$;
whence, providing $\lambda$ is a rational number, the angle $\varphi$ will
always be able to be determined algebraically from the angle $\psi$,
and consequently, since from the assumed angle $\varphi$ it will be
$u=\tang\psi=\sqrt{nnz^{2\lambda+2}-1}$,
everything will be able to be determined by this angle $\psi$,
because we will have
\[
nnz^{2\lambda+2}=1+\tang \psi^2=\frac{1}{\cos \psi^2}
\]
and hence $z=\sqrt[\lambda+1]{\frac{1}{n\cos\psi}}$;
then indeed, since $\varphi=\frac{\psi}{\lambda+1}$, the coordinates will be
$x=z\sin\frac{\psi}{\lambda+1}$ and $y=z\cos\frac{\psi}{\lambda+1}$;
all these values will thus be able to be exhibited algebraically.

\end{document}